\documentclass[a4paper,12pt]{article}
\usepackage{amsmath,amsthm,amssymb,fontenc,color,textcomp,amsfonts,graphicx,graphics}

\begin{document}

\begin{center}{\bf RIEMANNIAN GEODESICS - AN ILLUSTRATION FROM THE CALCULUS OF VARIATIONS}\end{center}
\vskip 1.0truecm

\begin{center}{ Uchechukwu Opara }\end{center}
\begin{center}{ African University of Science and Technology (AUST), Abuja, Nigeria \\ e-mail: ucmiop@yahoo.com . }\end{center}
\vskip 1.0truecm
\textbf{Abstract.} \ \ This paper sheds light on the essential characteristics of geodesics,  which frequently occur in considerations from motion in Euclidean space.  Focus is mainly on a method of obtaining them from the calculus of variations, and an explicit geodesic computation for a Riemannian hypersurface. \\ \\
\textbf{Keywords - } Riemannian Hypersurfaces, Geodesics, Calculus of Variations, Banach Path Spaces. \\ \\

\section {\bf Introduction.} \ \ The manifold of reference in this paper is the 3-sphere, sometimes referred to as the glome.  It is the focal point of the famous Poincar\'{e} conjecture, which since its proposal, has received continuous attention in the form of computational and analytic devices being crafted in an effort to study it.  One example of such devices is the \textit{Ricci Flow}, which was introduced into mainstream advanced mathematics in the mid 1980's.  A method for Ricci flow analysis is that of using \textit{minimal submanifolds} to probe the inner geometry of the manifolds in question. Understanding that geodesics are the basic class of minimal submanifolds, their involvement gives a relevant aspect into the study of a trending topic in modern Pseudo-Riemannian geometry. The higher dimensional minimal submanifold class of minimal hypersurfaces of 3-manifolds, which is more commonly consulted in Ricci Flow analysis, is better understood by knowledge of geodesics.  To buttress this fact, we may make reference to some of the copious pertinent publications in academic archives \cite{ac}, \cite{ae}, \cite{af}. The accuracy of the computational results of this paper, linked to a specific detailed choice of co-ordinate system, will yield much insight into the study of more general minimal submanifolds of the glome.\\

Generically, we may consider smooth, closed hypersurfaces embedded in real Euclidean space for a broad vantage point. In particular, for non-flat hypersurfaces with non-zero Gaussian curvature, it is evident that straight line segments cannot be admitted within the restrictions of their intrinsic structure.  Each geodesic curve must then be identified by way of anayltic devices.\\

Let $E$ be a real $n$-dimensional Euclidean space ($n \geq 3$) and $M \subset E$ be an ($n - 1$)-dimensional smooth connected submanifold, in other words, a hypersurface in $E$.  Let $c \subset M$ be a curve admitted by $M$ and parametrized by arclength or a constant multiple of it. Then $c$ is a geodesic of $M$ if and only if $\nabla _{\overset{\bf .}{c}} \overset{\bf .}{c} \equiv 0$, where the velocity vector field $\overset{\bf .}{c}$ along $c$ is simply the derivative of $c$ with respect to the parameter, and $\nabla$ is the covariant derivative.  The intrinsic differential operator $\nabla$ is given by the action; 
$$\nabla _X Y = D_X Y - \langle  D_X Y, \nu \rangle \nu $$
where $X$ and $Y$ are tangential vector fields to the hypersurface, $D$ is the directional derivative operator and $\nu$ is the outward unit normal vector field to $M$.\\

Clearly, $\overset{\bf .}{c}$ is a tangential vector field to $M$ and 
$$ \begin{array}{ccl}
\nabla _{\overset{\bf .}{c}} \overset{\bf .}{c} \equiv 0 & \Leftrightarrow & D _{\overset{\bf .}{c}} \overset{\bf .}{c} \equiv  \langle  D_{\overset{\bf .}{c}} \overset{\bf .}{c}, \nu \rangle \nu \\
{}& \Leftrightarrow & D _{\overset{\bf .}{c}} \overset{\bf .}{c} \equiv  \lVert  D_{\overset{\bf .}{c}} \overset{\bf .}{c} \rVert \nu
\end{array}$$
Observe that $c$ is also a submanifold with $\overset{\bf .}{c}$ as the spanning tangential basis vector field. As such,  $ D_{\overset{\bf .}{c}} \overset{\bf .}{c} = \overset{\bf ..}{c}$ and the above equivalence leads us to  $$\nu = \dfrac{ \overset{\bf ..}{c}}{\lVert  \overset{\bf ..}{c} \rVert} \ \ \mbox{whenever} \ \ \overset{\bf ..}{c} \neq 0.$$\\
We have $ \overset{\bf .}{c}$ in the direction of the principal unit tangent $T$ to $c$ and the principal unit normal $N$ to $c$ is obtained by the formula $$N = \dfrac{T'(s)}{\lVert T'(s) \rVert}$$ where $s$ is the arclength parameter.  But, recalling our chosen parameter, 
$$  \overset{\bf ..}{c}(ks) = k^2  \overset{\bf ..}{c}(s) = k^2  T'(s) $$ for any constant $k$, so that $$\dfrac{ \overset{\bf ..}{c}}{\lVert  \overset{\bf ..}{c} \rVert} = N.$$   
As such, the outward unit normal $\nu$ to $M$ and the principal unit normal $N$ to a geodesic $c \subset M$ coincide at each point of $c$.  We take this as an equivalent of the earlier stated geodesic characterization $\nabla _{\overset{\bf .}{c}} \overset{\bf .}{c} \equiv 0$. \\

With the analogy to motion of rigid bodies along $M$, taking the time parameter as a constant multiple of the arclength parameter, we get that for a geodesic $c \subset M$, either the acceleration vector  $\overset{\bf ..}{c}$ is null or the tangential component to $M$ is null, since $\nu$ is parallel to $ \overset{\bf ..}{c}$.  Of course, $ \overset{\bf ..}{c}$ can only be null in the event of motion along a straight line segment.\\ \\ \\
\section {\bf Obtaining Geodesic Solutions.} \ \ We thus far have the means of identifying the geodesics of hypersurfaces.  We now discuss a method of computing the geodesics for a given hypersurface.  Typically, we have two settings for such computations, namely;
\begin{description}
\item[i.)]compute a geodesic curve connecting two given points on the hypersurface,
\item[ii.)] compute a geodesic curve given a starting point and tangential direction on the hypersurface.
\end{description}

Although we are guaranteed the existence of solutions in either case due to the Hopf-Rinow theorem, in the first setting we typically do not have uniqueness of solution.  Let $p$ and $q$ be two given points on $M$, so that we denote by $\Omega(M;p,q)$ the set of all smooth curves on $M$ from $p$ to $q$.  This set, referred to as the path space of $M$ from $p$ to $q$, is itself an infinite-dimensional manifold and this is where all relevant optimization techniques are initiated.  This is to say the geodesic computation formulas are derived by differentiation in infinite dimensional (Banach) path spaces.   The curves in  $\Omega(M;p,q)$ are regarded as points of the infinite-dimensional manifold.\\

Let $\Gamma \in  \Omega(M;p,q)$ be such that $L(\Gamma) = \underset{c \in  \Omega(M;p,q)}{min} L(c)$, where $L(c)$ is the length of the curve $c$.  It is well-known that $\Gamma$ is a geodesic of $M$ but other possible curves in the path space  $\Omega(M;p,q)$ which are of greater length and still satisfy the criteria of being smooth geodesic curves are local minimizers of the arclength function in  $\Omega(M;p,q)$.  This is to say that for each geodesic  $ \Gamma ' \in \Omega(M;p,q)$ there exists a relatively open set $U$ in $\Omega(M;p,q)$ containing $\Gamma '$ for which $L(\Gamma ') = \underset{c \in U}{min} L(c).$ When the geodesic solution is unique, then it gives us the minimizer of the arclength function within the path space being considered.  For some clarity, an illustration of the path space for the second setting is hereby provided.  Let $$S := \{ \varphi (I) \subseteq \mathbb{R}^n : \varphi \in C^2(\overline{I}, \mathbb{R}^n)\} $$ and define $$\varphi _1(I) + \varphi _2(I) := \{\varphi _1(t) + \varphi _2(t) : t \in I \}$$ for a bounded subinterval $I$ of the reals, so that $S$ is a linear space.  As a hypersurface; $M$ may be given by the formula $g_M \ ^{-1}\{ k\}$ where $k \in \mathbb{R}$ is constant and $g_M$ is a smooth functional on $\mathbb{R}^n$.  In this event, the path space  $\Omega(M;p,{\bf v}) \subset S$  equals the set $$\{ \varphi (I) \in S : \varphi (inf I) = p \ ; \ \dfrac{d\varphi }{dt}|_{inf I} =\left\| \dfrac{d\varphi }{dt}|_{inf I} \right\| .{\bf v} \ ; \ g_M(\varphi (t)) = k \ \forall t \in I\}$$ and this is precisely the set of all smooth curves on $M$ starting from $p \in M$ in the unit direction ${\bf{v}}_p$ from the tangent space $  T_p M$.  $\Omega(M;p,\bf{v})$ is also infinite dimensional, but each geodesic in $\Omega(M;p,\bf{v})$ is either an extension or restriction of a unique $\Gamma \in \Omega(M;p,\bf{v})$.  Moreover, if $M$ can be represented by a periodic parametrization then we have a unique periodic extension for $\Gamma$ on $M$. We denote lengths of curves $c \in \Omega(M;p,\bf{v})$ by $L_{I,t}(c)$ where $t$ is a particular parameter running through the interval $I$.\\

Several constants of integration emerge in the course of geodesic computations, and because of guaranteed uniqueness of these constants with respect to a given co-ordinate system for the second case  $\Omega(M;p,\bf{v})$, this setting is better amenable than the previously discussed path space  $\Omega(M;p,q)$. \\

The arclength function in the n-dimensional space $E$ is given by the indefinite integral $$\int \sqrt{\sum_{i = 1}^n {\pi _i}^2} \ ,$$ where $\{\pi _i\}_{i=1}^n$ is the set of canonical projection maps on $E$.  Taking $E$ to be $\mathbb{R}^n$ with the Euclidean co-ordinate system $(x_1 , x_2 , \cdots , x_n)$ then the projection maps $\pi _i$ are written as $dx _i$ for $i = 1, 2, \cdots , n$ and consequently, the arclength function is given by  $$s = \int \sqrt{\sum_{i = 1}^n {dx _i}^2} \ .$$  Expressing each $x_i$ in terms of a common parameter $t$ running through an interval $I \subseteq \mathbb{R}$ (which we can always do for a smooth curve in Euclidean space), we get the arclength function := $s_{I,t}$ as a definite integral $$s_{I,t} = \int_I  \sqrt{\sum_{i = 1}^n \left({\dfrac{dx_i}{dt}}\right)^2} dt \ .$$ We have $s_{I,t}$ as a functional map on the path space $\Omega(M; p, \bf{v})$; 
$$\begin{array}{rclc}
s_{I,t} : & \Omega(M; p, \bf{v}) & \rightarrow & \mathbb{R}\\
{} & c & \mapsto & L_{I,t}(c),
\end{array}$$
where $L_{I,t}(c)$ retains its previous notation.\\

The arclength function $s_{I,t}$  on the considered path spaces is convex and continuously differentiable, but the lack of reflexivity of the path spaces thwarts arguments about uniqueness of geodesic solutions by conventional optimization theorems.  However, in light of the uniqueness of the geodesic (and thus also the minimizer) in our current setting, any relevant optimization device employed here will yield the unique solution.\\

For instance, the canonical Euler - Lagrange equations  are derived by differentiation in the infinite-dimensional path spaces, and they are popularly used for this type of problem.  We have a functional of the type $$s_{I,t} = \int _I \Lambda (t, u_1 (t), u_2 (t), \cdots , u_{n-1}(t), \overset{\bf .}{u_1}, \overset{\bf .}{u_2}, \cdots, \overset{\bf .}{u}_{n-1})dt$$ taking $M$ to be parametrized by the co-ordinate system $$f : U \subseteq \mathbb{R}^{n-1} \rightarrow M \subset \mathbb{R}^n \ ; \ (u_1, u_2, \cdots , u_{n-1}) \mapsto (x_1 , x_2 , \cdots , x_n) \ ,$$ where $U \subseteq \mathbb{R}^{n-1}$ is a cube, or a convex domain for which $\overset{\circ}{U} \neq \emptyset$.  By the canonical Euler-Lagrange equations, all minimizing points of $s_{I,t}$ must satisfy $$\Lambda_{u_i} = \dfrac{d}{dt}(\Lambda_{\overset{\bf .}{u_i}})$$ for $1 \leq i \leq n-1$.  Here, $\overset{\bf .}{u_i}$ denotes $\dfrac{du_i}{dt}$ and $Y_x$ denotes partial derivation of $Y$ with respect to the variable connoted by the subscript $x$.  We obtain the transformation of $s_{I,t}$ from its initial form $ \displaystyle \int_I  \sqrt{\sum_{i = 1}^n \left({\dfrac{dx_i}{dt}}\right)^2} dt $ by way of the Riemannian structure $\{ g_{ij} \}_{i,j = 1}^{n-1}$ obtained from the co-ordinate system $f$ where $g_{ij} = \langle \dfrac{\partial f}{\partial u_i}, \dfrac {\partial f}{\partial u_j} \rangle$.  The function $s_{I,t}$ then becomes $$ \int _I \sqrt{\sum_{i,j = 1}^{n-1} g_{ij} \overset{\bf .}{u_i} \overset{\bf .}{u_j}}dt$$ which is the suitable form to apply the Euler - Lagrange equations.  \\

We can derive the weak formulation for the Euler-Lagrange equations by considering path spaces one more time.  This time we work in $$S^* := \{ \varphi (I) \subseteq \mathbb{R}^{n-1} : \varphi \in C^1(\overline{I}, \mathbb{R}^{n-1})\}$$ with the norm $$||\varphi|| = ||\varphi||_{\infty}+||\overset{\bf .}{\varphi} ||_{\infty}$$ where $$||x||_{\infty} = \underset{t \in \overline{I}}{sup}||x(t)||_2$$ giving us a Banach structure on $(S^*,||.||)$.  Let $p = f(a)$ and $Df_p(\psi) := f'(a ; \psi )= {\bf v}$.  Hence, we are minimizing the functional $s_{I,t}$ over the convex subset $\Omega(f^{-1}(M); a,\psi) \subseteq S^*$ which is an important advantage of computing in this setting instead of in the non-convex $\Omega (M , p, {\bf v}) \subset S$.  Assuming that $\overline{\varphi}$ is a local minimizer of $s_{I,t}$ in  $\Omega(f^{-1}(M); a,\psi) \subseteq S^*$ then $$\exists r > 0 : s_{I,t}(\overline{\varphi}) \leq s_{I,t} (\varphi ) \ \forall \varphi \in  \Omega(f^{-1}(M); a,\psi ) \cap B(\overline{\varphi} , r) \ .$$  In this setting, for any $\zeta \in S^*$ satisfying $\zeta (inf I) = \zeta (sup I) = \dfrac{d \zeta}{dt}|_{inf I} = 0$ we have $\overline{\varphi} + \tau \zeta \in \Omega(f^{-1}(M); a,\psi)$ where $\tau$ runs through an open interval containing zero.  Also, there exists some real number $\delta > 0$ such that for all $\tau \in (-\delta , \delta)$, we have $\overline{\varphi} + \tau \zeta \in B(\overline{\varphi} , r)$.  Define $\gamma (\tau) := s_{I,t}(\overline{\varphi} + \tau \zeta )$ so that $\gamma (0) \leq \gamma(\tau) \ \forall \tau \in (-\delta , \delta)$.  In other words, $0$ is a minimizer of $\gamma$ in $(-\delta , \delta)$ which means $\gamma '(0) = 0 \ \Rightarrow \ s_{I,t}'(\overline{\varphi}).\zeta = 0 .$   \\ \\ \\ $\begin{array}{lcl} s_{I,t}'(\overline{\varphi}).\zeta &=& \underset{\alpha \rightarrow 0}{lim}\left(\dfrac{s_{I,t}(\overline{\varphi} + \alpha \zeta) - s_{I,t}(\overline{\varphi})}{\alpha}\right) \\ {}&=& \displaystyle{\int}_I \ \underset{\alpha \rightarrow 0}{lim}\dfrac{\Lambda (t, \overline{\varphi}(t) + \alpha \zeta (t) , \overset{\bf .}{\overline{\varphi}}(t) + \alpha \overset{\bf .}{\zeta}(t)) - \Lambda (t , \overline{\varphi}(t) , \overset{\bf .}{\overline{\varphi}}(t))}{\alpha}dt \end{array}$ \\ \\
The quantity under the above integral is identified as the following directional derivative in $\mathbb{R} \times \mathbb{R}^{n-1} \times \mathbb{R}^{n-1}$ \\ \\ $\Lambda '(t , \overline{\varphi}(t) , \overset{\bf .}{\overline{\varphi}}(t));(0, \zeta (t) , \overset{\bf .}{\zeta} (t)) \\ = \Lambda _u(t, \overline{\varphi}(t) ,  \overset{\bf .}{\overline{\varphi}}(t)).\zeta (t) + \Lambda _{\overset{\bf .}{u}}(t, \overline{\varphi}(t) , \overset{\bf .}{\overline{\varphi}}(t)).\overset{\bf .}{\zeta} (t)$ \\  \\ where $u = (u_1 , u_2 , \cdots , u_{n-1})$ and $\overset{\bf .}{u} = (\overset{\bf .}{u_1}, \overset{\bf .}{u_2}, \cdots, \overset{\bf .}{u}_{n-1})$.  We justify passing the limit into the above integral by way of uniform convergence of the integrand at differentiable points of $\Lambda$  in the limit as $\alpha \rightarrow 0$.  Resultantly, the weak formulation of the Euler-Lagrange equations is $$\int _I \left( \Lambda _u(t, \overline{\varphi}(t) ,  \overset{\bf .}{\overline{\varphi}}(t)).\zeta (t) + \Lambda _{\overset{\bf .}{u}}(t, \overline{\varphi}(t) ,  \overset{\bf .}{\overline{\varphi}}(t)).\overset{\bf .}{\zeta} (t) \right) dt = 0$$  for all $\zeta \in S^*$ satisfying $\zeta (inf I) = \zeta (sup I) = \dfrac{d \zeta}{dt}|_{inf I} = 0$.  Integrating by parts and assuming in addition that $\overline{\varphi}$ is of class $C^2$, we have $$\int _I \left( \Lambda _u(t, \overline{\varphi}(t) ,  \overset{\bf .}{\overline{\varphi}}(t)) - \dfrac{d}{dt} \Lambda _{\overset{\bf .}{u}}(t, \overline{\varphi}(t) ,  \overset{\bf .}{\overline{\varphi}}(t)) \right).\zeta dt = 0 \ .$$  By the basic lemma of the calculus of variations, the first vector in the integrand; $\Lambda _u(t, \overline{\varphi}(t) ,  \overset{\bf .}{\overline{\varphi}}(t)) - \dfrac{d}{dt} \Lambda _{\overset{\bf .}{u}}(t, \overline{\varphi}(t) ,  \overset{\bf .}{\overline{\varphi}}(t)) \in \mathbb{R}^{n-1}$ vanishes giving us $$\Lambda_{u_i} = \dfrac{d}{dt}(\Lambda_{\overset{\bf .}{u_i}})$$ for $1 \leq i \leq n-1$. Hence the computed curve lies in $\mathbb{R}^{n-1}$ while the geodesic is its image under $f$.\\ 
\indent Note that geodesics are intrinsic to their manifolds, and so are invariant under which co-ordinate systems are used to derive them.  Nevertheless, a suitable choice of co-ordinate system could significantly ease a computation.  Different manifolds and their co-ordinate systems present their own challenges but the main theoretical properties are tantamount in all applications. \\

In motion along manifolds, minimal geodesics are the `lazy' curves which minimize both the arclength function and energy functionals.  For the path space $\Omega (M; p, {\bf v})$, a kinetic energy functional is given by 
$$ \begin{array} {rclc}
U_{I,t} :& \Omega(M; p, {\bf v}) & \rightarrow & \mathbb{R}\\
{} & c & \mapsto & \displaystyle \int _I \lVert \overset{\bf .}{c}(t) \rVert^2 dt
\end{array}$$
and another is given by $$c \longmapsto \int _I  {\kappa}_{c}^2 ds \ ,$$
where $\kappa _c$ is the curvature of $c$ in terms of the arclength parameter.\\

 In essence, for a particle at rest on a smooth Riemannian hypersurface, when a linear tangential force is applied to it so that it undergoes motion on the hypersurface without slipping, it is translated along the geodesic path.  In vector calculus, the phenomenon of geodesics is central to several themes, including parallel transports and Jacobi fields. \\ 
\section {\bf Computation of Geodesic Equation for Unit 3-Sphere in $ \mathbb{R}^4$ .}
With all the given background theoretical information, we now proceed to demonstrate the computation of a hyperspherical geodesic.  According to a concise argument by Jost and Li-Jost\cite{af} (p.42), all such curves must be portions of great circles of the hypersphere.  We will operate using the hyperspherical co-ordinate system, which yields geodesic equations easier to solve in comparison to most other known co-ordinate systems.  The minimal global parametrization of the unit sphere $S^{n-1} \subset \mathbb{R}^n$ by the hyperspherical co-ordinate system for $n \geq 3$ is given by \\

${f}_{(n-1)}:{[\frac{-\pi}{2},\frac{\pi}{2}]}^{n-2} \times [0,2\pi] \longrightarrow {S}^{n-1} \ \ ;
({u}_{1},{u}_{2},\cdots,{u}_{n-1}) \longmapsto ({x}_{1},{x}_{2},\cdots,{x}_{n}); \\ \\
\indent {x}_{1} = \substack{n-1\\ \prod\\ j=1} \cos {u}_{j} \ \ ; \ \  \ {x}_{i} = \substack{n - i\\ \prod\\ j=1} \cos {u}_{j}\sin {u}_{n-i+1}$ \ [$2 \leq i \leq n-1], \\
\indent {x}_{n} = \sin {u}_{1}$ .\\

For the case of interest of least dimension ($S^2$), we have a multiple parametrization\\

$\overline{f_2}:[\frac{-\pi}{2},\frac{\pi}{2}] \times \mathbb{R} \longrightarrow \mathbb{R}^3;\\
\indent \overline{f_{2}}^1(u_1, u_2) := x = \cos u_1 \cos u_2 \\ 
\indent \overline{f_{2}}^2(u_1, u_2) := y = \cos u_1 \sin u_2\\
\indent \overline{f_{2}}^3(u_1, u_2) := z = \sin u_1 $ .\\

\noindent  Pre-images of geodesic solutions computed in $[\frac{-\pi}{2}, \frac{\pi}{2}] \times \mathbb{R}$ are known to have a parametrization $\left( u_1 , \arctan \left[\dfrac{\sin u_1}{\sqrt{\dfrac{\cos ^2 u_1}{\gamma ^2} -1}} \right] + \delta \right)$ for real constants of integration $\delta$ and $\gamma,  | \gamma | < 1$.  This solution holds only for $|u_1| < \arccos(| \gamma |)$.  Singularities occur at $u_1 = \pm \arccos(| \gamma |)$, where we respectively take $u_2 = \pm \frac{\pi}{2} + \delta$ and for $|u_1| > \arccos(| \gamma |)$, we have no real solutions.  This is due to the fact that the solution curve does not cut across the parallels ${u_1}^{*}$ of $S^2$ for which $|{u_1}^{*}| > \arccos(| \gamma |)$.  We extend the solution into another branch $$\left \{ \left( u_1 , \pi - \arctan \left[\dfrac{\sin u_1}{\sqrt{\dfrac{\cos ^2 u_1}{\gamma ^2} -1}} \right] + \delta \right) \right \}_{| u_1 | < \arccos(| \gamma |)}$$ to get one full cycle of the geodesic curve, considering the minimal domain of $\overline{f_2}$.  For either branch, the image computed under $\overline{f_2}$ yields one half of a great circle of ($S^2$).  The integration constant $\gamma$ equals 0 if and only if the solution curve is a meridian and the special case of $| \gamma | = 1$ occurs if and only if the solution is the equator. \\

We will now proceed with the computation of a geodesic equation for the unit hypersphere $S^3 \subset \mathbb{R}^4$.  It is given once by the global co-ordinate system,
$$\begin{array}{rlcl}
f_3 : & [\frac{-\pi}{2}, \frac{\pi}{2}]^2 \times [0, 2\pi] & \longrightarrow & S^3 \subset \mathbb{R}^4 \\
{}& (u_1 , u_2 , u_3) & \longmapsto & f(u_1 , u_2 , u_3) = (x_1 , x_2 , x_3 , x_4)
\end{array}$$
$f_3 (u_1 , u_2 , u_3) = (\cos u_1 \cos u_2 \cos u_3 , \cos u_1 \cos u_2 \sin u_3 , \cos u_1 \sin u_2 , \sin u_1)$ .\\

\noindent $S^3$ is also covered multiple times by extending $f_3$ to $\overline{f_3} := f$ with a domain $[\frac{-\pi}{2},\frac{\pi}{2}] \times \mathbb{R}^2$ and the same parametrization in the co-domain $\mathbb{R}^4$.  In both cases, using the notation $g_{ij} = \langle \dfrac{\partial f}{\partial u_i}, \dfrac {\partial f}{\partial u_j} \rangle$, we get\\

\noindent $g_{11} = 1\\
g_{22} = \cos ^2 u_1 \\
g_{33} = \cos ^2 u_1 . \cos ^2 u_2 \\
g_{ij} = 0, \ \ \  \mbox{for} \  i\neq j.$\\

Setting the common parameter to be $u_1$, we get an arclength function on the path space $\Omega(S^3; p, {\bf v})$ for ${\bf v}_p \in T_p S^3$; 
$$ s_{I, u_1} = \int _I \sqrt{1 + \cos ^2 u_1 \left( \dfrac{du_2}{du_1} \right) ^2 + \cos ^2 u_1 \cos ^2 u_2 \left( \dfrac{du_3}{du_1} \right) ^2} du_1 := \int _I \Lambda du_1$$
for an appropriate interval $I \subseteq [\frac{-\pi}{2}, \frac{\pi}{2}]$. In determining the mininimizer $c$ of $s_{I, u_1}$, our choice of parameter ideally ensures that $||c'(u_1)|| \neq 0$ on $I$ and it leaves us to set up the following Euler - Lagrange equations (1) and (2)- -\\

For the first Euler - Lagrange equation,
$$\begin{array}{rcl}
\Lambda _{u_3} & =  &\dfrac{d}{du_1}( \Lambda _{\overset{\bf .}{u_3}})\\
\Leftrightarrow  \left( \dfrac{du_3}{du_1} \right) ^2 & = & \dfrac{\gamma ^2 \left( \sec ^2 u_1 + (\frac{du_2}{du_1})^2 \right)}{\cos ^2 u_1 \cos ^4 u_2 - \gamma ^2 \cos ^2 u_2} \ \ \ \cdots (1)
\end{array}$$

for a real constant $\gamma, | \gamma | < 1$.\\

For the second Euler - Lagrange equation,
$$ \Lambda _{u_2}  =  \dfrac{d}{du_1}( \Lambda _{\overset{\bf .}{u_2}}) $$ 
$$ \Leftrightarrow \dfrac{-1}{\Lambda} ( \cos ^2 u_1 \cos  u_2 \sin u_2 )(\overset{\bf .}{u_3})^2 =  \dfrac{d}{du_1}\left( \dfrac{ \overset{\bf .}{u_2} \cos ^2 u_1}{\Lambda} \right)$$
 $$=  \dfrac{\Lambda (-2 \overset{\bf .}{u_2} \cos u_1 \sin u_1 + \overset{\bf ..}{u_2} \cos ^2 u_1) - \overset{\bf .}{u_2}(\frac{d\Lambda}{du_1})\cos ^2 u_1}{\Lambda ^2}  \ \ \ \cdots (2). $$

Consistent with previous notation, $\overset{\bf .}{Y} = \dfrac{dY}{du_1}$ and $\overset{\bf ..}{Y} = \dfrac{d^2Y}{du_1 ^2}$. \\

After grinding out with lengthy arithmetic and substituting $(1)$ in $(2)$, we are able to obtain a differential equation in exactly 2 variables: \\ \\
  $\overset{ \bf .}{u_2}\sin u_1 \cos u_2 (\gamma ^2 - 2\cos ^2 u_1 \cos ^2 u_2 ) + \overset{ \bf ..}{u_2}\cos u_1 \cos u_2 ( \cos ^2 u_1 \cos ^2 u_2 - \gamma ^2) +  \gamma ^2 \sec u_1 \sin u_2 + \gamma ^2 (\overset{\bf .}{u_2})^2 \cos u_1 \sin u_2 - (\overset{\bf .}{u_2})^3 \cos ^4 u_1 \sin u_1 \cos ^3 u_2 = 0 \ \ \ \cdots (3) $ \\ \\
It is useful to say how to obtain $\gamma$ by way of simple vector calculus.  Already having ${\bf v}_p \in T_p S^3$ and the gradient of $f$ at $p, Df_p$, we obtain the pre-image of the position vector ${\bf v}_p$ under $Df$ which we denote $ {\bf\psi}$ so that $Df_p(\psi) = {\bf v}_p$.  Now, having the vector $\psi := (\psi _1, \psi _2, \psi _3)$ and the pre-image of the point $p$ under $f$,  we substitute $du_3(\psi) = \psi _3, du_2(\psi) = \psi _2$ and $du_1(\psi) = \psi _1$ in $(1)$ to realize the value of $\gamma$. ($\{ du_i \}_{i=1}^3$  is the set of canonical projection maps on the domain $\mathbb{R}^3$.)\\

Upon inspection, $u_2 \equiv 0$ is a solution of the non-linear differential equation (3).  This solution occurs if and only if the $u_2$ co-ordinate at the initial position is 0 and the pullback of the initial velocity vector ${\bf v}_p$ has a null second component, which is to say $\psi _2 = 0$. To find the general solution for (3), we re-arrange the equation as follows --  \\ \\
$[ 1 + (\frac{du_2}{du_1})^2 \cos ^2 u_1 ] \ [\gamma ^2 \sec ^2 u_1 \tan u_2 - (\frac{du_2}{du_1})\tan u_1 \cos ^2 u_1 \cos^2 u_2] = $\\
$[ \gamma ^2 - \cos ^2 u_1 \cos ^2 u_2] \ [\frac{d^2 u_2}{du_1 ^2} - (\frac{du_2}{du_1})\tan u_1]$ \\ 
\vspace{0.5 cm}
\indent The equation then becomes easier to decompose, \\
$\begin{array}{lccl}
{}& \dfrac{ \overset{\bf ..}{u_2} - \overset{\bf .}{u_2} \tan u_1}{ 1 + (\overset{\bf .}{u_2})^2 \cos ^2 u_1} & = & \dfrac{ \gamma ^2 \sec ^2 u_1 \tan u_2 - \overset{\bf .}{u_2} \tan u_1 \cos ^2 u_1 \cos ^2 u_2 }{\gamma ^2 - \cos ^2 u_1 \cos ^2 u_2 }\\ \\
\Leftrightarrow & \dfrac{ \frac{d}{du_1}(1 + (\overset{\bf .}{u_2})^2 \cos ^2 u_1)}{2 \overset{\bf .}{u_2} \cos ^2 u_1 [ 1 + (\overset{\bf .}{u_2})^2 \cos ^2 u_1]} & = &  \dfrac{ \gamma ^2 \sec ^2 u_1 \tan u_2 - \overset{\bf .}{u_2} \tan u_1 \cos ^2 u_1 \cos ^2 u_2 }{\gamma ^2 - \cos ^2 u_1 \cos ^2 u_2 }\\ \\
\Leftrightarrow & \dfrac{ \frac{d}{du_1}(1 + (\overset{\bf .}{u_2})^2 \cos ^2 u_1)}{ 1 + (\overset{\bf .}{u_2})^2 \cos ^2 u_1 } & = &  \dfrac{ 2\gamma ^2 \overset{\bf .}{u_2}\tan  u_2 \sec ^2 u_2 - 2 (\overset{\bf .}{u_2})^2 \tan u_1 \cos ^4 u_1  }{\gamma ^2 \sec ^2 u_2 - \cos ^2 u_1 }
\end{array}$\\ \\ \\ 
\indent $=  \dfrac{ 2\gamma ^2 \overset{\bf .}{u_2}\tan  u_2 \sec ^2 u_2 + 2\cos u_1 \sin u_1 - 2 \cos u_1 \sin u_1  [ 1 + (\overset{\bf .}{u_2})^2 \cos ^2 u_1] }{\gamma ^2 \sec ^2 u_2 - \cos ^2 u_1 }$\\ \\
\indent $=  \dfrac{ \frac{d}{du_1}(\gamma ^2 \sec ^2 u_2 - \cos ^2 u_1) }{\gamma ^2 \sec ^2 u_2 - \cos ^2 u_1 } -  \dfrac{  2 \cos u_1 \sin u_1  [ 1 + (\overset{\bf .}{u_2})^2 \cos ^2 u_1] }{\gamma ^2 \sec ^2 u_2 - \cos ^2 u_1 }$\\ \\ \\
Integrating both sides with respect to $u_1$ , we get
$$ ln[ 1 + (\overset{\bf .}{u_2})^2 \cos ^2 u_1 ] = ln[ \gamma ^2 \sec ^2 u_2 - \cos ^2 u_1 ] - \int \dfrac{  2 \cos u_1 \sin u_1  [ 1 + (\overset{\bf .}{u_2})^2 \cos ^2 u_1] }{\gamma ^2 \sec ^2 u_2 - \cos ^2 u_1 } du_1$$
$$ \Leftrightarrow ln \left( \dfrac{ 1 + (\overset{\bf .}{u_2})^2 \cos ^2 u_1 }{\gamma ^2 \sec ^2 u_2 - \cos ^2 u_1 }\right) = -\int 2 \cos u_1 sin u_1 \left( \dfrac{ 1 + (\overset{\bf .}{u_2})^2 \cos ^2 u_1 }{\gamma ^2 \sec ^2 u_2 - \cos ^2 u_1 }\right)du_1$$ \\
Hence, we solve for $U$ in the equation; \\ \\
$\begin{array}{lccl}
{} &\displaystyle \int U(t).V(t) dt & = & ln(U(t)) \ \ \ \ [t = u_1] \\ \\
\Leftrightarrow & U.V & = & \left(\dfrac{1}{U} \right) \dfrac{dU}{dt} \\ \\
\Leftrightarrow & \dfrac{\frac{dU}{dt}}{U^2} & = & V\\ \\
\Leftrightarrow & -\dfrac{1}{U} & = &\displaystyle \int V dt
\end{array}$ \\ \\ \\
Swapping back $U$ with \ $\dfrac{ 1 + (\overset{\bf .}{u_2})^2 \cos ^2 u_1 }{\gamma ^2 \sec ^2 u_2 - \cos ^2 u_1 }$ \ and $V$ with $-2 \cos u_1 \sin u_1$, \\ \\ we get (3) as an equivalent separable differential equation - - \\
$\begin{array}{lccl}
{} & \dfrac{\gamma ^2 \sec ^2 u_2 - \cos ^2 u_1 }{ 1 + (\overset{\bf .}{u_2})^2 \cos ^2 u_1 } & = & -\cos ^2 u_1 + k \\ \\
\Leftrightarrow & \gamma ^2 \sec ^2 u_2 & = & - (\overset{\bf .}{u_2})^2 \cos ^4 u_1 + k + k (\overset{\bf .}{u_2})^2 \cos ^2 u_1 \\ \\
\Leftrightarrow & \gamma ^2 \sec ^2 u_2 - k & = & \left( \dfrac{du_2}{du_1} \right)^2 [ k\cos ^2 u_1 - \cos ^4 u_1 ] \ \ \ \cdots (3')
\end{array}$ \\ \\
The constant of integration $k$ above is obtained in a similar way as $\gamma$. \\  \\ (*)Apologetically, we will fix $k$ here to equal $2 \gamma ^2 - \gamma ^4$ in order to shorten the computation and avoid carrying along too many constants of integration.  With this substitution, we have as a non-trivial solution for $(3')$;\\
$u_2 = \arctan \left[\dfrac{\sin u_1}{\sqrt{\dfrac{\cos ^2 u_1}{\gamma ^2} -1}} \right] = \arcsin \left[ \dfrac{\tan u_1}{\sqrt{\dfrac{1}{\gamma ^2} - 1}} \right]$ \ \ \ \ for $|u_1| < \arccos(\sqrt{2 \gamma ^2 - \gamma ^4})$\\
$\Rightarrow \dfrac{du_2}{du_1} = \dfrac{\gamma \sec u_1}{\sqrt{ \cos ^2 u_1 - \gamma ^2}}$ \ \ \ and \\
\indent $\cos u_2 = \left( \dfrac{\frac{1}{\gamma ^2} - \sec ^2 u_1}{ \frac{1}{\gamma ^2} - 1} \right)^{\frac{1}{2}}$.\\ \\
Henceforth, we will write $t$ instead of $u_1$.  From equation (1), we have
$$ \left(\dfrac{du_3}{dt} \right)^2 = \left( \sec ^2 t + \dfrac{\gamma ^2 \sec ^2 t}{\cos ^2 t - \gamma ^2} \right)\left(\dfrac{\gamma ^2}{\cos ^2 t \left(\frac{1 - \gamma ^2 \sec ^2 t}{1 - \gamma ^2}\right)^2 - \frac{\gamma ^2 - \gamma ^4 \sec ^2 t}{1 - \gamma ^2 } } \right)$$
$$ \Rightarrow \left| \dfrac{du_3}{dt} \right| = \dfrac{|\gamma|(1 - \gamma ^2)}{\sqrt{(\cos t - \gamma ^2 \sec t)^2 [ \cos ^2 t - 2 \gamma ^2 + \gamma ^4]}} \ \ \ \cdots (4)$$ 
A solution for (4) when $|t| < \arccos(\sqrt{2 \gamma ^2 - \gamma ^4})$ is
$$ u_3 = \arcsin \left( \dfrac{ \tan t}{\sqrt{\frac{1}{\gamma ^2} - 1}\sqrt{1 - \gamma ^2 \sec ^2 t}} \right) + \beta $$
The constant $\beta$ is determined by obtaining the pre-image of the starting point $p$ under $f$ and substituting in the solution above.  The above solution given for $u_3$ runs through an interval of length $\pi$.  Given the minimal domain of $f$, we extend this solution into another branch, namely - - 
$$\overline{u_3} = \pi - \arcsin \left( \dfrac{ \tan t}{\sqrt{\frac{1}{\gamma ^2} - 1}\sqrt{1 - \gamma ^2 \sec ^2 t}} \right) + \beta $$
so that the third argument of the domain runs through an interval of length $2 \pi$, hence completing one full cycle of the solution curve.  Let us set $z = u_3 - \beta$ (resp. $z = \overline{u_3} - \beta$) so that we have;\\ 
$\begin{array}{rl}
\cos t \cos u_2 \cos u_3 = & \cos t \cos u_2 ( \cos z \cos \beta - \sin z \sin \beta)\\  
=&\pm \cos t  \left( \dfrac{\frac{1}{\gamma ^2} - \sec ^2 t}{ \frac{1}{\gamma ^2} - 1} \right)^{\frac{1}{2}}  \left( \dfrac{\frac{1}{\gamma ^2} + (\gamma ^2 -2) \sec ^2 t}{(\frac{1}{\gamma ^2} - 1)(1 - \gamma ^2 \sec ^2 t )}\right) ^{\frac{1}{2}} \cos \beta - \\ {} & \cos t  \left(\dfrac{\frac{1}{\gamma ^2} - \sec ^2 t}{ \frac{1}{\gamma ^2} - 1} \right)^{\frac{1}{2}} (\sin z) (\sin \beta) \\
=& \dfrac{\pm \cos t \left(\frac{1}{\gamma ^2} + (\gamma ^2 -2) \sec ^2 t \right )^{\frac{1}{2}} \cos \beta}{\gamma ( \frac{1}{\gamma ^2} -1)} - \dfrac{ \sin t \sin \beta}{\gamma (\frac{1}{\gamma ^2} - 1)}\\ \\
\cos t \cos u_2 \sin u_3 = & \cos t \cos u_2 ( \cos z \sin \beta + \sin z \cos \beta)\\ 
= &\dfrac{\pm \cos t \left(\frac{1}{\gamma ^2} + (\gamma ^2 -2) \sec ^2 t \right )^{\frac{1}{2}} \sin \beta}{\gamma ( \frac{1}{\gamma ^2} -1)} +  \dfrac{ \sin t \cos \beta}{\gamma (\frac{1}{\gamma ^2} - 1)}\\
\cos t \sin u_2 = & \dfrac{\sin t}{\sqrt{\frac{1}{\gamma ^2}-1}}
\end{array}$\\ \\
Setting the parameters ${\varphi _1|}_{S^3} = \sin t$ and ${\varphi _2|}_{S^3} = \cos t \left(\frac{1}{\gamma ^2} + (\gamma ^2 -2) \sec ^2 t \right )^{\frac{1}{2}}$, we see that the solution curve lies on the 2-plane in $\mathbb{R}^4$ parametrized by 
$$ \left( \dfrac{(\cos \beta) \varphi _2 - (\sin \beta) \varphi _1}{\gamma ( \frac{1}{\gamma ^2} -1)} \ , \  \dfrac{(\sin \beta) \varphi _2 + (\cos \beta) \varphi _1}{\gamma ( \frac{1}{\gamma ^2} -1)} \ , \  \dfrac{\varphi _1}{\sqrt{\frac{1}{\gamma ^2}-1}} \ , \  \varphi _1 \right) .$$  
Since this plane passes through the origin in $\mathbb{R}^4$, its intersection with $S^3$, which is precisely the locus of our solution curve, is a great circle of the hypersphere as expected. The curve does not intersect sections of the hypersurface for which $|u_1| > \arccos(\sqrt{2 \gamma ^2 - \gamma ^4})$ and it is left as an exercise for the reader to verify that $$L_{I , u_1}(c) := \int _{|u_1| < \arccos(\sqrt{2 \gamma ^2 - \gamma ^4})} \lVert Df (\overset{\bf .}{c}) \rVert du_1 = \pi$$ holds for either solution branch $c$ in the domain of $f$.\\ \\
The graph below depicts the pre-image of our solution curve under $f$ in $\mathbb{R}^3$ for $\beta = 0$ and $\gamma = \frac{\sqrt{2}}{2}$.\\
\includegraphics[scale=0.4]{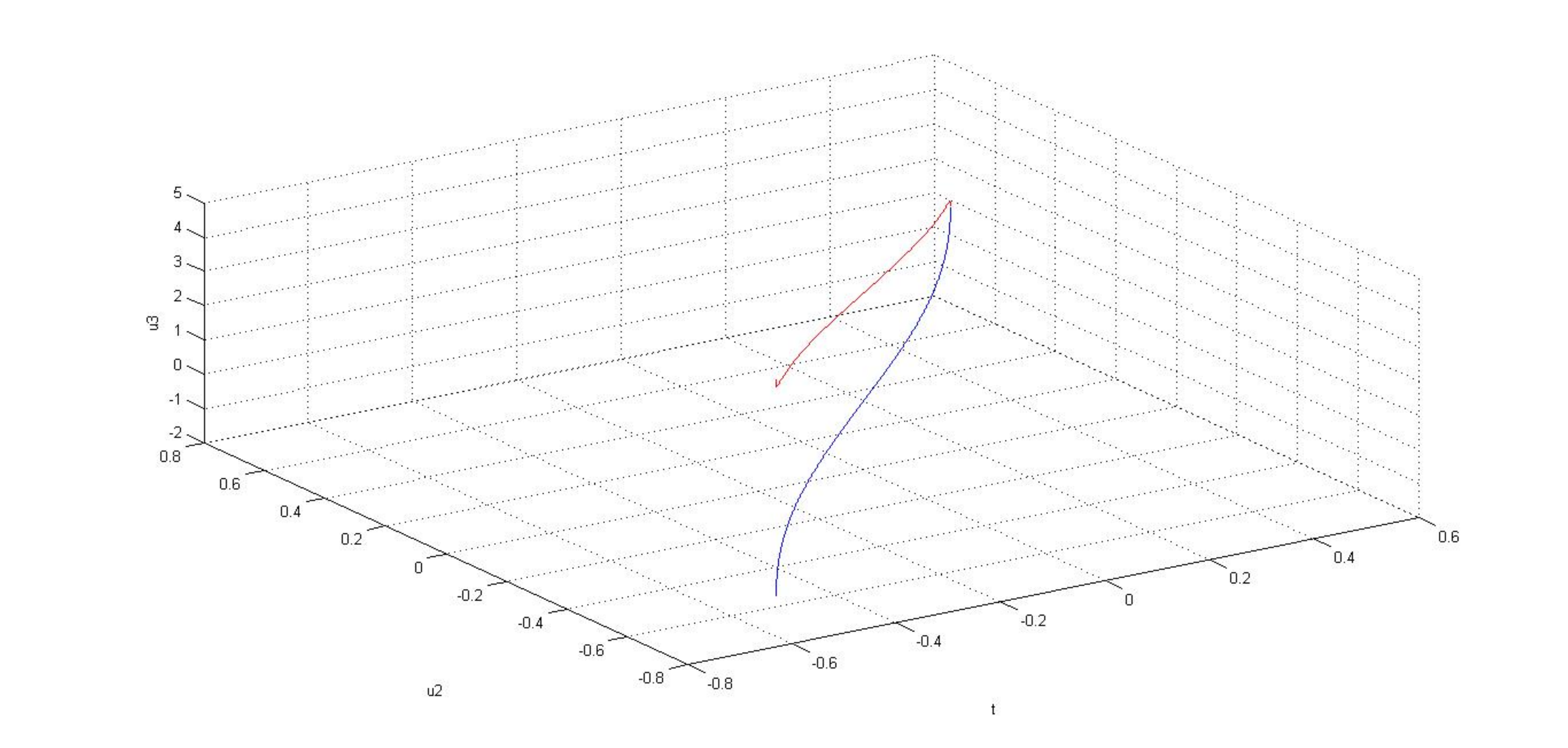} \\
\indent Due to the action taken earlier in this section (*) dropping another constant of integration, the solution we have corresponds to only a single integral curve of equation $(3')$.  All geodesics of $S^3$ correspond to the integral curves of $(3')$, which depend on the constant $\gamma$ to be determined beforehand.  With these observations, we have the extrapolation of the 2-spherical geodesic equation into one higher dimension.  Although somewhat tedious in leaping from the second dimension which is of greatest interest, this should not be regarded as a foolhardy venture as all results are obtainable in terms of elementary functions and hold potential for beneficial analysis.

 \end{document}